\documentclass[psamsfonts]{amsart}

\usepackage{amssymb,amsfonts}
\usepackage[all,arc]{xy}
\usepackage{enumerate}
\usepackage{mathrsfs}

\newtheorem{theorem}{Theorem}[section]
\newtheorem{corollary}[theorem]{Corollary}
\newtheorem{proposition}[theorem]{Proposition}
\newtheorem{lemma}[theorem]{Lemma}

\theoremstyle{definition}
\newtheorem{definition}[theorem]{Definition}

\theoremstyle{remark}
\newtheorem{remark}[theorem]{Remark}

\DeclareMathOperator{\ord}{ord}

\DeclareMathOperator{\cont}{cont}

\makeatletter

\begin{document}
	\title{Gauss' lemma for polynomials over semidomains}
	
	\author[P. Nasehpour]{Peyman Nasehpour}
	
	\address{Department of Engineering Science \\
		Golpayegan University of Technology \\ Golpayegan \\ Iran}
	\email{nasehpour@gut.ac.ir, nasehpour@gmail.com}
	
	

	\begin{abstract}
	In this paper, we generalize Gauss' lemma for polynomials over subtractive factorial semidomains.
	\end{abstract}
	
	\maketitle
	
	\section{Introduction}

In this paper, by a semiring, we understand an algebraic structure, consisting of a nonempty set $S$ with two operations of addition and multiplication such that the following conditions are satisfied:

\begin{enumerate}
	\item $(S,+)$ is a commutative monoid with an identity element $0$;
	\item $(S,\cdot)$ is a commutative monoid with an identity element $1 \not= 0$;
	\item Multiplication distributes over addition, i.e., $a\cdot (b+c) = a \cdot b + a \cdot c$ for all $a,b,c \in S$;
	\item The element $0$ is the absorbing element of the multiplication, i.e., $s \cdot 0=0$ for all $s\in S$.
\end{enumerate}

Since the language for semirings is not completely standardized \cite{Glazek2002}, we need to introduce some other concepts. A nonempty subset $I$ of a semiring $S$ is said to be an ideal of $S$, if $a+b \in I$ for all $a,b \in I$ and $sa \in I$ for all $s \in S$ and $a \in I$ \cite{Bourne1951}. An ideal $I$ of a semiring $S$ is called a proper ideal of the semiring $S$, if $I \neq S$. A proper ideal $P$ of a semiring $S$ is called a prime ideal of $S$, if $ab\in P$ implies either $a\in P$ or $b\in P$. Finally, let us recall that an ideal $I$ of a semiring $S$ is subtractive if $a+b\in S$ and $a\in S$ imply that $b\in S$ for all $a,b\in S$ \cite{Golan1999(b)}.

A semiring $S$ is a semidomain if $ab=ac$ with $a\neq 0$ will cause $b=c$, for all $a,b,c\in S$. Similar to the concept of field of fractions in ring theory, one can define the semifield of fractions $F(S)$ of the semidomain $S$ \cite[p. 22]{Golan1999(a)}. An ideal $I$ of a semiring $S$ is called principal if $I = \{sa: s\in S\}$ for some $a\in S$. The ideal $I = \{sa: s\in S\}$ is denoted by $(a)$. Finally, if $S$ is a semiring, for $a,b \in S$, it is written $a \mid b$ and said that ``$a$ divides $b$'', if $b = sa$ for some $s\in S$. This is equivalent to say that $(b) \subseteq (a)$. Also, it is said that $a$ and $b$ are associates if $a=ub$ for some unit and note that if $S$ is a semidomain, then this is equivalent to say that $(a) = (b)$. A nonzero, nonunit element $s$ of a semiring $S$ is said to be irreducible if $s = s_1 s_2$ for some $s_1, s_2 \in S$, then either $s_1$ or $s_2$ is a unit. This is equivalent to say that $(s)$ is maximal among proper principal ideals of $S$. An element $p\in S-\{1\}$ is said to be a prime element, if the principal ideal $(p)$ is a prime ideal of $S$, which is equivalent to say if $p \mid ab$, then either $p \mid a$ or $p \mid b$ \cite{Nasehpour2018}.

A semidomain $S$ is called factorial (also unique factorization) if the following conditions are satisfied:

\begin{enumerate}
	
	\item Each irreducible element of $S$ is a prime element of $S$.
	
	\item Any nonzero and nonunit element of $S$ is a product of irreducible elements of $S$.
	
\end{enumerate}
	
In Section \ref{sec:content}, we introduce the concept of the order of an element in a factorial (unique factorization) semidomain and the traditional version of the concept of the content of a polynomial over such semirings, all inspired by the Lang's approach to these concepts in his textbook (refer to Chapter IV $\S2$ in \cite{Lang2002}). 

In Definition \ref{orderatp}, we define that if $S$ is a factorial semidomain, $F$ is its semifield of fractions, $a$ is a nonzero element of $F$, and $p$ is a prime element of $S$, then the order of $a$ at $p$ is the integer $v$, denoted by $\ord_p(a) =v$, such that $a= p^v \cdot b$ and $p$ does not divide the numerator or denominator of $b\in F$. If $a = 0$, we define its order at $p$ to be $\infty$. With the help of this concept, we define the concept of the ``order at $p$'' for polynomials in one variable. For the polynomial \[f=a_0 + a_1 X + \cdots + a_n X^n \in F[X],\] if $f = 0$, we define $\ord_p (f) = \infty$ and $f \neq 0$, we define $\ord_p(f)$ to be \[ \ord_p(f) = \min_{i \in \Lambda} \ord_p a_i,\] where $\Lambda = \{i: 0 \leq i\leq n, a_i \neq 0\}$. Then, we define the content of $f$, denoted by $\cont(f)$, to be the following product \[\prod p^{\ord_p(f)}, \] being taken over all $p$ such that $\ord_p(f) \neq 0$, or any multiple of this product by a unit of $S$.

It is straightforward to see that for each $f \in F[X]$, there is a polynomial $f_1 \in F[X]$ such that $f = c f_1$, $c=\cont(f)$, and $\cont (f_1) =1$ (see Definition \ref{traditionalcontentdef} and Proposition \ref{traditionalcontentprop}).

In Section \ref{sec:gauss}, we prove a semiring version of Gauss' lemma for polynomial semirings (see Theorem \ref{Gausslemma}) in this sense that the function $\cont: F[X] \rightarrow F$ is multiplicative, i.e. for all $f,g \in F[X]$, we have \[\cont(fg) = \cont(f) \cont(g).\]

Our general references for semiring theory are the books \cite{Golan1999(a),Golan1999(b),Golan2003}.
	
\section{The Concepts of Order and Content}\label{sec:content}
	
Let $S$ be a factorial semidomain and $F$ its semifield of fractions. Let $a\in F$ and $a\neq 0$. If $p$ is a prime element of $S$, then $a$ can be written as $a= p^v \cdot b$ such that $v$ is an integer number and $p$ does not divide the numerator or denominator of $b\in F$. Since $S$ is factorial, $v$ is uniquely determined by $a$. Based on this simple argument, we give the following definition:
	
	\begin{definition}
		
		\label{orderatp}
		
		Let $S$ be a factorial semidomain and $F$ its semifield of fractions. Let $a$ be a nonzero element of $F$. If $p$ is a prime element of $S$, then we define the order of $a$ at $p$ the integer $v$, denoted by $\ord_p(a) =v$, such that $a= p^v \cdot b$ and $p$ does not divide the numerator or denominator of $b\in F$. If $a = 0$, we define its order at $p$ to be $\infty$.
	\end{definition}

\begin{definition}
	Let $S$ and $T$ be two semirings. We say that a function $f: S-\{0\} \rightarrow T$ has logarithmic property if $f(xy) = f(x)+f(y)$ for all nonzero $x,y \in S$.
\end{definition}
	
	 The proof of the following statement is straightforward:
	
	\begin{proposition}
		Let $S$ be a factorial semidomain, $F$ its semifield of fractions, and $p$ a prime element of $S$. Then the following statements hold:
		
		\begin{enumerate}
			\item The function $\ord_p: F-\{0\} \rightarrow \mathbb Z$ has logarithmic property, i.e. \[\ord_p(xy) = \ord_p(x) + \ord_p(y),\] where $x,y\in F$ and $xy \neq 0$.
			
			\item For all nonzero $a\in S$, we have $\ord_p(a) =1$ if and only if $p\nmid a$. 
		\end{enumerate}
	\end{proposition}
	
	\begin{remark}
		Let $S$ be a factorial semidomain, $F$ its semifield of fractions, and $p$ a prime element of $S$. In fact, $\ord_p$ is a discrete valuation map on $F$ (see Definition 1.1 and Example 3.2 in \cite{NV2018}). 
	\end{remark}
	
	Now, we define the concept of the ``order at $p$'' for polynomials in one variable:    
	
	\begin{definition}
		
		\label{traditionalcontentdef}
		
		Let $S$ be a factorial semidomain, $F$ its semifield of fractions, and $p$ a prime element of $S$. Let \[f=a_0 + a_1 X + \cdots + a_n X^n \in F[X]\] be a polynomial in one variable. 
		
		\begin{enumerate}
			\item If $f = 0$, we define $\ord_p (f) = \infty$. If $f \neq 0$, we define $\ord_p(f)$ to be \[ \ord_p(f) = \min_{i \in \Lambda} \ord_p a_i,\] where $\Lambda = \{i: 0 \leq i\leq n, a_i \neq 0\}$. 
			
			\item If $v=\ord_p(f)$, we call $up^v$ a $p$-content for $f$, for any unit $u$ of $S$.
			
			\item We define the content of $f$, denoted by $\cont(f)$, to be the following product \[\prod p^{\ord_p(f)}, \] being taken over all $p$ such that $\ord_p(f) \neq 0$, or any multiple of this product by a unit of $S$.
		\end{enumerate}
		
	\end{definition}
	
	Let $S$ be a semiring. A greatest common divisor (abbreviated as gcd) of a set $A \subseteq S$, which has at least one nonzero element, is a nonzero element $d$, if $d \mid a$ for any $a\in A$ and if $d^{\prime} \mid a$ for any $a\in A$, then $d^{\prime} \mid d$. It is clear that this is equivalent to say that $(d)$ is the minimal element of all principal ideals containing the ideal generated by the set $A$. A greatest common divisor of the set $A$, which is not necessarily unique, is denoted by $\gcd(A)$ \cite{Nasehpour2018}. The proof of the following statement is straightforward:
	
	\begin{proposition}
		
		\label{traditionalcontentprop}
		
		Let $S$ be a factorial semidomain, $F$ its semifield of fractions, and $f \in F[X]$. Then the following statements hold:
		
		\begin{enumerate}
			\item The content of $f$ is well defined up to multiplication by a unit of $S$.
			
			\item The function $\cont: F[X] \rightarrow F$ is homogeneous, i.e. for each nonzero $b\in F$, we have $\cont (bf) = b \cont (f)$ \cite{NasehpourGMJ}.
			
			\item There is a polynomial $f_1 \in F[X]$ such that $f = c f_1$, $c=\cont(f)$, and $\cont (f_1) =1$. In particular, all coefficients of $f_1$ lie in $S$ and their gcd is 1. 
		\end{enumerate}
	\end{proposition}

	\section{Gauss' Lemma for Polynomials over Subtractive Factorial Semidomains}\label{sec:gauss}
	
	\begin{lemma}
		\label{dividesubtractive}
		Let $S$ be a subtractive semiring and $d$, $a$, and $b$ be arbitrary elements of $S$. If $d\mid a$ and $d\mid a+b$, then $d\mid b$.
		
		\begin{proof} 
			Let $d\mid a$ and $d\mid a+b$. So, $(a) \subseteq (d)$ and $(a+b) \subseteq (d)$. Therefore, $(a+b,a) \subseteq (d)$. Since $S$ is subtractive and $a$ and $a+b$ are elements of the ideal $(a+b,a)$, we have $b\in (a+b,a)$. This implies that $(b) \subseteq (d)$, which is equivalent to say that $d\mid b$ and this finishes the proof.
		\end{proof}
		
	\end{lemma}
	
	Now, we have enough tools to give a semiring version of Gauss' Lemma:
	
	\begin{theorem}[Gauss' Lemma]
		
		\label{Gausslemma}
		
		Let $S$ be a subtractive factorial semidomain, and let $F$ be its semifield of fractions. Let $f,g \in F[X]$ be polynomials in one variable. Then the function $\cont: F[X] \rightarrow F$ is multiplicative, i.e. for all $f,g \in F[X]$, we have \[\cont(fg) = \cont(f) \cont(g).\]
		
		\begin{proof}
			By Proposition \ref{traditionalcontentprop}, we can write $f=cf_1$ and $g=dg_1$ such that $c =\cont(f)$, $d=\cont(g)$, $\cont(f_1)=1$, and $\cont(g_1)=1$. Therefore, it is sufficient to prove that if $\cont(f)=1$ and $\cont(g)=1$, then $\cont(fg) =1$, and for this, it suffices to prove that for each prime $p$, $\ord_p(fg)=0$. Now, let \[f=a_n X^n + \cdots + a_0, \text{with~} a_n \neq 0, \text{and,} \] \[g = b_m X^m + \cdots + b_0, \text{with~} b_m \neq 0\] be polynomials of content 1. Let $p$ be a prime of $S$. In order to show that $\ord_p(fg)=0$, it will suffice to prove that $p$ does not divide all coefficients of $fg$. Let $r$ be the largest integer such that $0 \leq r \leq n$, $a_r \neq 0$, and $p$ does not divide $a_r$. Similarly, let $b_s$ be the coefficient of $g$ farthest to the left, $b_s \neq 0$, such that $p$ does not divide $b_s$. Consider the coefficient of $X^{r+s}$ in $fg$. This coefficient is equal to \[c_{r+s} = a_r b_s + a_{r+1} b_{s-1} + \cdots + a_{r-1} b_{s+1} + \cdots \] and $p \nmid a_r b_s$. However, $p$ divides every other non-zero term in this sum since in each term there will be some coefficient $a_i$ to the left of $a_r$ or some coefficient $b_j$ to the left of $b_s$. Now, since every ideal of $S$ is subtractive, by Lemma \ref{dividesubtractive}, $p$ can not divide $c_{r+s}$, and this completes the proof of the lemma.\end{proof}
	\end{theorem}

	\begin{definition}
		Let $S$ be a factorial semidomain and $F$ its semifield of fractions. We define a polynomial of content 1 to be a primitive polynomial.
	\end{definition}

	\begin{corollary}
		Let $S$ be a factorial semidomain. If $f,g\in S[X]$ are primitive, then $fg$ is also primitive.
	\end{corollary}

	\begin{remark}
		The modern definition for the concept of the content of an element of a semialgebra has been given and discussed in \cite{Nasehpour2016}. If $S$ is a factorial semidomain and $f\in S[X]$ is a polynomial, then the content of $f$ in the modern form, denoted by $c(f)$, is defined to be the ideal generated by the coefficients of $f$ and in such a case, it is easy to see that $c(f)$ is the principal ideal generated by $\cont (f)$. An $S$-semialgebra $B$ is called Gaussian if $c(fg)=c(f)c(g)$, for all $f,g \in B$. It is now clear that by Gauss' Lemma (Theorem \ref{Gausslemma}) if $S$ is a factorial semidomain, then $S[X]$ is a Gaussian $S$-semialgebra.
	\end{remark}

	\subsection*{Acknowledgments} The author is supported by the Department of Engineering Science at the Golpayegan University of Technology and his special thanks go to the Department for providing all necessary facilities available to him for successfully conducting this research.

\end{document}